\documentclass[a4paper,11pt]{amsart}

\usepackage{amssymb}
\usepackage{mathrsfs}
\usepackage{amsmath,amssymb,amsthm,latexsym,amscd,mathrsfs}
\usepackage{indentfirst}
\usepackage{stmaryrd}
\usepackage{graphicx}
\usepackage{extarrows}
\usepackage{multirow}
\usepackage{latexsym}
\usepackage{amsfonts}
\usepackage{color}
\usepackage{pictexwd,dcpic}
\usepackage{graphicx}
\usepackage{psfrag}
\usepackage{hyperref}
\usepackage{comment}
\usepackage{enumerate}

\usepackage{cite}

\numberwithin{equation}{section} \theoremstyle{plain}
\newtheorem{thm}{Theorem}[section]

\newtheorem{prop}[thm]{Proposition}

\newtheorem{conj}[thm]{Conjecture}

\newtheorem{rem}[thm]{Remark}

\newtheorem{ques}[thm]{Question}

\newtheorem{ack}{Acknowledgements}   
 \makeatletter

\def\<{\langle}
\def\>{\rangle}
\def\({\left(}
\def\){\right)}
\def\[{\left[}
\def\]{\right]}

\DeclareMathOperator{\diag}{diag}
\DeclareMathOperator{\rank}{rank}
\DeclareMathOperator{\Tr}{Tr}

\makeatother

\title{A survey on the DDVV-type inequalities}

\author[J.Q. Ge]{Jianquan Ge}
\address{School of Mathematical Sciences, Laboratory of Mathematics and Complex Systems, Beijing Normal University, Beijing 100875, P.R. CHINA.}
\email{jqge@bnu.edu.cn}

\author[F.G. Li]{Fagui Li}
\address{Beijing International Center for Mathematical Research, Peking University,
Beijing 100871, P.R. CHINA.}
\email{faguili@bicmr.pku.edu.cn}

\author[Z. Z. Tang]{Zizhou Tang$^{\dag}$}
\address{Chern Institute of Mathematics \& LPMC, Nankai University, Tianjin 300071, P. R. China}
\email{zztang@nankai.edu.cn}

\author[Y. Zhou]{Yi Zhou}
\address{Beijing International Center for Mathematical Research, Peking University,
Beijing 100871, P.R. CHINA.}
\email{yizhou@bicmr.pku.edu.cn}

\date{}
\thanks {$^{\dag}$ the corresponding author}
\thanks{J. Q. Ge is partially supported by NSFC (No. 12171037) and the Fundamental Research Funds for the Central Universities.}
\thanks{F. G. Li is partially supported by  NSFC (No. 12171037, 12271040) and China Postdoctoral Science Foundation (No. 2022M720261).}
\thanks{Z. Z. Tang is partially supported by NSFC (Nos. 11931007, 11871282), Nankai Zhide Foundation, Tianjin Outstanding Talents Foundation.}
\thanks{Y. Zhou is partially supported by  NSFC (No. 12171037, 12271040), China Postdoctoral Science Foundation (No. BX20230018)
and National Key R$\&$D Program of China 2020YFA0712800.}

\begin{document}

\begin{abstract}
In this paper, we give a survey on the history and recent developments on the DDVV-type inequalities.
\end{abstract}
\maketitle

\tableofcontents
\section{Introduction of DDVV-type inequality}
A DDVV-type inequality is an estimate of the form
\begin{equation}\label{DDVVtypeineq}
\sum^m_{r,s=1}\|[B_r,B_s]\|^2\leq c\(\sum^m_{r=1}\|B_r\|^2\)^2,
\end{equation}
considered for certain type of $n\times n$ matrices $B_1, \ldots, B_m$,
where $[A,B]:=AB-BA$ denotes the commutator, $\|B\|^2:=\operatorname{tr}(BB^*)$ denotes the squared Frobenius norm ($B^*$ be the conjugate transpose)
and $c$ is a nonnegative constant.
The DDVV-type inequality originates from the normal scalar curvature conjecture
(which is also called the DDVV conjecture) in submanifold geometry.

In 1999, De Smet-Dillen-Verstraelen-Vrancken \cite{DDVV99}
proposed the normal scalar curvature conjecture:
\begin{conj}[DDVV Conjecture \cite{DDVV99}]
Let $M^n$ be an immersed submanifold of a real space form with constant sectional curvature $\kappa$. Then
\begin{equation}\label{geoDDVVineq}
\rho+\rho^{\bot}\leq |H|^2+\kappa,
\end{equation}
where $\rho$ denotes the normalized scalar curvature,
$\rho^{\bot}$ denotes the normalized normal scalar curvature
and $H$ denotes the normalized mean curvature vector field.
\end{conj}
Several years later, Dillen-Fastenakels-Veken \cite{DFV07} pointed out that
the geometric inequality (\ref{geoDDVVineq}) is ture if the algebraic inequality (\ref{DDVVtypeineq})
holds for arbitrary $n\times n$ real symmetric matrices $B_1, \ldots, B_m$
with the universal constant $c=1$ (which is called the DDVV inequality).
After some partial results were obtained \cite{CL08, DFV07a, DFV07, DHTV04, Lu07},
this inequality on real symmetric matrices was finally proved
by Lu \cite{Lu11} and Ge-Tang \cite{GT08} independently and differently.
Here, we state their theorem as follows.
Define an action of $K(n, m):=\operatorname{O}(n)\times\operatorname{O}(m)$
on a family of matrices $(B_1,\cdots,B_m)$ by
$$(P,R)\cdot(B_1,\cdots,B_m):=(P^*B_1P,\cdots,P^*B_mP)\cdot R.$$

\begin{thm}[DDVV inequality \cite{GT08, Lu11}]\label{algDDVV}\label{DDVVsymthm}
Let $B_1,\cdots,B_m$ be arbitrary $n\times n$ real symmetric matrices $(m,n\geq2)$.
Then
\begin{equation}\label{algDDVVineq}
\sum^m_{r,s=1}\|\[B_r,B_s\]\|^2\leq \(\sum^m_{r=1}\|B_r\|^2\)^2.
\end{equation}
The equality holds if and only if there exists a $(P,R)\in K(n, m)$ such that
$$(P,R)\cdot(B_1,\cdots,B_m)=
(\diag(H_1,0), \diag(H_2,0), 0, \cdots, 0),$$
where for some $\lambda\geq0$,
$$H_1:=
\begin{pmatrix}
\lambda & 0 \\
0 & -\lambda \\
\end{pmatrix},\
H_2:=
\begin{pmatrix}
0 & \lambda \\
\lambda & 0 \\
\end{pmatrix}.$$
\end{thm}

A natural question is to determine the optimal universal constant $c$ in (\ref{DDVVtypeineq})
when the matrices $B_1, \ldots, B_m$ are in other regarded classes.
Using the same approach \cite{GT08, GT11},
Ge \cite{Ge14} proved the DDVV-type inequality for real skew-symmetric matrices,
and applied it to Yang-Mills fields in Riemannian submersion geometry.
Subsequently, Ge-Xu-You-Zhou \cite{GXYZ17} extended the DDVV-type inequalities from real symmetric and skew-symmetric matrices to Hermitian and skew-Hermitian matrices.
Further, by decomposing complex matrices into Hermitian matrices and skew-Hermitian matrices,
Ge-Li-Zhou \cite{GLZ21} generalized the DDVV-type inequalities for complex (skew-)symmetric matrices
and general complex (real) matrices.
According to the same philosophy,
Ge-Li-Zhou \cite{GLZ21} also obtained the DDVV-type inequality for general quaternionic matrices.

In 2005, B\"ottcher-Wenzel \cite{BW05} raised the so-called BW conjecture:
\begin{conj}[BW conjecture \cite{BW05}]
Let $X, Y$ be arbitrary $n\times n$ real matrices, then
\begin{equation}\label{BWineq}
\|\[X,Y\]\|^2\leq 2\|X\|^2\|Y\|^2.
\end{equation}
\end{conj}
Since $4\|X\|^2\|Y\|^2\leq (\|X\|^2+\|Y\|^2)^2$,
the above BW inequality $($\ref{BWineq}$)$ implies the DDVV inequality for the case $m=2$.
B\"{o}ttcher-Wenzel \cite{BW05} and L\`{a}szl\`{o} \cite{L07} first proved
the special case of $n=2$ and $n=3$, respectively.
And then the complete proofs were given by B\"{o}ttcher-Wenzel \cite{BW08}, Vong-Jin \cite{VJ08}, Audenaert \cite{AKMR09} and Lu \cite{Lu11, Lu12} in various ways.
B\"{o}ttcher-Wenzel \cite{BW08} also extended the BW inequality from real matrices to complex matrices,
and Cheng-Vong-Wenzel \cite{CVW10} obtained the characterization of the equality.
By introducing the Kronecker product, Ge-Li-Lu-Zhou \cite{GLLZ20} obtained new proofs
of the complex BW inequality and its equality condition.
For the convenience of readers, we restate the BW inequality and the equality condition in the following.
\begin{thm}[BW inequality \cite{LW17}]
Let $X, Y$ be arbitrary real (or complex) matrices, then $$\|\[X,Y\]\|^2\leq 2\|X\|^2\|Y\|^2, $$
where the equality holds  if and only if there exists a unitary matrix $U$ such that $X = U(X_0\oplus O)U^*$ and $Y = U(Y_0 \oplus O)U^*$
with a maximal pair $(X_0, Y_0)$ of $2\times2$ matrices. Here, $O$ denotes the zero matrix in corresponding order, two matrices $X,Y \in  \mathbb{C}^{2\times2}$ form a maximal pair if and only if $\operatorname{tr}X=\operatorname{tr}Y=0$ and $\left\langle X, Y\right\rangle  = 0$.
\end{thm}
Also, analogues with, e.g., Schatten norms
were investigated (cf. \cite{CFL13, CL17, W10, WKMRA10}, etc).
 In addition,
Ge-Li-Zhou \cite{GLZ21} generalized the BW inequality from complex matrices to quaternionic matrices.
Some other generalizations of the BW-type inequalities were also obtained by
Fong-Cheng-Lok\cite{FCL10}, Wenzel\cite{W10}, Wenzel-Audenaert \cite{WKMRA10},
Cheng-Fong-Lei \cite{CFL13}, Cheng-Liang \cite{CL17}, Cheng-Akintoye-Jiao \cite{CAJ21}, L\`{a}szl\`{o}\cite{L22} and Chru\'{s}ci\'{n}ski-Kimura-Ohno-Singal\cite{CKOS22, CKOS23}, etc. For more details, we recommend \cite{CJV15} and \cite{LW17} for a comprehensive overview on the developments.

We summarize the optimal constant $c$ for DDVV-type inequalities and BW-type inequalities mentioned above in the following Table \ref{table1} and Table \ref{table2}.
\begin{table}[h]
\caption{The optimal constant $c$ of DDVV-type inequalities ($m\geq3$)} \label{table1}
\centering
\begin{tabular}{|c|c|c|}
\hline
c &real &complex\\
\hline
symmetric & $1$ &$1$\\
\hline
skew-symmetric &$\frac13(n=3), \frac23(n\geq4)$ &$\frac13(n=3), \frac23(n\geq4)$ \\
\hline
Hermitian & --  &$\frac43$\\
\hline
skew-Hermitian & -- &$\frac43$\\
\hline
general &$\frac43$ &$\frac43$\\
\hline
\end{tabular}
\end{table}

\begin{table}[h!]
\caption{The optimal constant $c$ of DDVV-type and BW-type inequalities}\label{table2}
\centering
\begin{tabular}{|c|c|c|c|}
\hline
c &\ \ \ \ \ \ real\ \ \ \ \ \ &\ \ \ complex\ \ \ &quaternionic\\
\hline
DDVV$(m\geq3)$ &$\frac43$ &$\frac43$ &$\frac83$\\
\hline
DDVV$(m=2)$ &1 &1 &2\\
\hline
BW &2 &2 &4\\
\hline
\end{tabular}
\end{table}

\section{The origin and applications in geometry}
\subsection{DDVV inequality in submanifold geometry}
Suppose $M^n$ is an immersed submanifold of a real space form $N^{n+m}(\kappa)$
with constant sectional curvature $\kappa$.
Let $R$ be the Riemannian curvature tensor of $M$, let
$R^{\bot}$ be the curvature tensor of the normal connection,
and let $\operatorname{II}$ be the second fundamental form.
For an arbitrary point $p\in M$,
let $\{e_1,\cdots,e_n\}$ and $\{\xi_1,\cdots,\xi_m\}$ be orthonormal bases of
$T_pM$ and $T_p^{\bot}M$, respectively.
Then the normalized scalar curvature is
$$\rho=\frac{2}{n(n-1)}\sum_{1=i<j}^n\langle R(e_i, e_j)e_j, e_i\rangle,$$
the normalized normal scalar curvature is
$$\rho^{\bot}=\frac{2}{n(n-1)}\Bigg(\sum_{1=i<j}^n\sum_{1=r<s}^m\langle R^{\bot}(e_i, e_j)\xi_r, \xi_s\rangle^2\Bigg)^{\frac12}=\frac{2}{n(n-1)}|R^{\bot}|$$
and the normalized mean curvature vector field is
$H=\frac1n\sum_{i=1}^{n}\operatorname{II}(e_i, e_i)$.
For each $1\leq r\leq m$,
let $A_r$ be the matrix correspond to the shape operator in direction $\xi_r$
with respect to the basis $\{e_1,\cdots,e_n\}$, and let $$B_r:=A_r-\langle H, \xi_r\rangle I_n.$$
On the one hand, by the Gauss equation, we have
$$-n(n-1)(\rho-\kappa)=|\operatorname{II}|^2-n^2|H|^2
=\sum^m_{r=1}\|B_r\|^2-n(n-1)|H|^2,$$ and thus
\begin{equation}\label{DDVVRHS}
|H|^2-\rho+\kappa=\frac1{n(n-1)}\sum^m_{r=1}\|B_r\|^2.
\end{equation}
On the other hand, by the Ricci equation, we have
\begin{equation}\label{DDVVLHS}
\rho^{\bot}=\frac1{n(n-1)}\Bigg(\sum^m_{r,s=1}\|[A_r,A_s]\|^2\Bigg)^{\frac12}
=\frac1{n(n-1)}\Bigg(\sum^m_{r,s=1}\|[B_r,B_s]\|^2\Bigg)^{\frac12}.
\end{equation}
It follows from (\ref{DDVVRHS}) and (\ref{DDVVLHS}) that
the geometric DDVV inequality (\ref{geoDDVVineq}) can be derived from
the algebraic DDVV inequality (\ref{algDDVVineq}).
Therefore, Theorem \ref{algDDVV} implies the following normal scalar curvature inequality:

\begin{thm}[Ge-Tang\cite{GT08}, Lu\cite{Lu11}]\label{DDVV}
Let $M^n$ be an immersed submanifold of a real space form $N^{n+m}(\kappa)$. Then
$$\rho+\rho^{\bot}\leq |H|^2+\kappa.$$
The equality holds at some point $p\in M$ if and only if there exist an orthonormal basis $\{e_1,\cdots,e_n\}$ of $T_pM$ and an orthonormal basis $\{\xi_1,\cdots,\xi_m\}$ of $T_p^{\bot}M$ such that
$$A_1=\lambda_1I_n+\mu\diag(H_1,0),\
A_2=\lambda_2I_n+\mu\diag(H_2,0),\
A_3=\lambda_3I_n$$ and $A_{\xi_r}=0$ for $r>3$,
where $\mu, \lambda_1, \lambda_2, \lambda_3$ are real numbers.
\end{thm}

There are many generalized normal scalar curvature inequalities under different geometric assumptions
such as special submanifolds in statistical manifolds, complex space forms or Sasakian space forms
(cf. \cite{ACDV20, AMM17, AO18, BUS19, CC11, LLV22, M14, M16, M17, R17, WX23, Zhan19}).


The DDVV inequality has many other important applications in submanifold geometry.
For example,
it was used in Gu-Xu's work \cite{GX12} on Yau rigidity theorem for minimal submanifolds in spheres.
The case of two symmetric matrices is a core technology in the proof of
the following integral inequality by Chern-do Carmo-Kobayashi \cite{CdK}.
To simplify notations, we denote by $S$ the squared length of the second fundamental form.

\begin{thm}[Chern-do Carmo-Kobayashi\cite{CdK}, Simons\cite{Simons}]\label{simonsineq}
Let $M^n$ be a closed, minimal, immersed submanifold of a real space form $N^{n+m}(\kappa)$. Then
$$\int_M\left[\left(2-\frac1m\right)S-n\kappa\right]S\ dV_M\geq0.$$
\end{thm}

The above integral inequality started a series of explorations on the gap phenomena and pinching results
for the second fundamental form (cf. \cite{Law69, NS69, Simons, CdK, Yau, Shen89, WS98, LL92, CX93}).
Finally, we have the following pinching theorem:

\begin{thm}[Chern-do Carmo-Kobayashi\cite{CdK}, Lawson\cite{Law69}]\label{hypsurpin}
Let $M^n$ be a minimal hypersurface in the unit sphere $\mathbf{S}^{n+1}$. 
If $0\leq S\leq n$, then $M$ is either totally geodesic or is one of the Clifford tori
$$M_{k,n-k}=\mathbf{S}^k\left(\sqrt{\frac{k}{n}}\,\right)\times\mathbf{S}^{n-k}\left(\sqrt{\frac{n-k}{n}}\,\right).$$
\end{thm}

\begin{thm}[Chen-Xu\cite{CX93}, Li-Li\cite{LL92}]\label{submfdpin}
Let $M^n$ be a closed, minimal submanifold in the unit sphere $\mathbf{S}^{n+m}$, $m\geq 2$.
If $0\leq S\leq\frac23 n$, then $M$ is either totally geodesic or is a Veronese surface in $\mathbf{S}^{2+m}$.
\end{thm}

In \cite{Lu11}, Lu considered
the fundamental matrix on $M$ which is an $m\times m$ matrix-valued function defined as $A=(a_{rs})$,
where $a_{rs}=\langle A_r, A_s\rangle$.
Let $\lambda_1\geq\cdots\geq\lambda_m\geq0$ be the eigenvalues of the fundamental matrix.
Note that $S$ is the trace of the fundamental matrix, i.e., $S=\lambda_1+\cdots+\lambda_m$. Set $\lambda_2:=0$ if $m=1$.
The following result with pinching quantity $S+\lambda_2$ was proved
by using a generalized DDVV inequality (see Theorem \ref{thm Lu inequality}).

\begin{thm}[Lu\cite{Lu11}]\label{Lupin}
Let $M^n$ be a closed minimal submanifold in the unit sphere $\mathbf{S}^{n+m}$.
If $$0\leq S+\lambda_2\leq n,$$
then $M$ is totally geodesic, or is one of the Clifford tori $M_{k,n-k}$ $(1\leq k< n)$ in $\mathbf{S}^{n+m}$,
or is a Veronese surface in $\mathbf{S}^{2+m}$.
\end{thm}

Since $\lambda_2\leq \frac12S$, the above Theorem \ref{Lupin} extended the former rigidity results
(Theorem \ref{hypsurpin} and Theorem \ref{submfdpin}).
Furthermore, Leng-Xu \cite{LX18} generalized Lu's rigidity theorem to submanifolds with parallel mean curvature.

\subsection{Simons-type integral inequality in Riemannian submersion geometry}
In some sense, Riemannian submersions can be seen as the ``dual" of isometric immersions,
while real skew-symmetric matrices can also be seen as a kind of ``dual" of real symmetric matrices.
Fortunately, analogous to the case of real symmetric matrices,
the DDVV-type inequality for real skew-symmetric matrices can also deduce
a Simons-type integral inequality for Riemannian submersions.
In fact, Ge \cite{Ge14} applied it to give a Simons-type inequality for (generalized) Yang-Mills fields in Riemannian submersions geometry (dual to Simons inequality for minimal submanifolds of spheres in submanifold geometry \cite{Lu11}). The dual phenomenon between Yang-Mills fields and minimal submanifolds was initially investigated by Tian \cite{Ti00}.

Let $\pi: M^{n+m}\rightarrow B^n$ be a Riemannian submersion,
and let $D$ be the Levi-Civita connection on $M$.
Then the O'Neill integrability tensor $A$ is defined by
$$A_{X}Y:=\mathscr{H}D_{\mathscr{H}X}\mathscr{V}Y+\mathscr{V}D_{\mathscr{H}X}\mathscr{H}Y,$$
where $\mathscr{H}$ and $\mathscr{V}$ denote the projections from the tangent bundle $TM$ to the horizontal and vertical distribution, respectively. $A$ is essentially the curvature when the Riemannian submersion is a Euclidean vector bundle projection, and thus a similar equation $\delta^DA=0$ about $A$ is used to define a Yang-Mills horizontal distribution as Yang-Mills connections on  Euclidean vector bundles.
Since the $2$-tensor field $A$ is alternating on the horizontal distribution, locally it can be represented by $m$ skew-symmetric matrices of order $n$.
We first state the DDVV-type inequality for real skew-symmetric matrices:

\begin{thm}[Ge\cite{Ge14}]\label{skewDDVV}
Let $B_1,\cdots,B_m$ be $n\times n$ real skew-symmetric matrices.
\begin{enumerate}
\item If $n=3$, then
\begin{equation*}
\sum^m_{r,s=1}\|\[B_r,B_s\]\|^2\leq \frac13\(\sum^m_{r=1}\|B_r\|^2\)^2.
\end{equation*}
The equality holds if and only if there exists a $(P,R)\in K(n, m)$ such that
$$(P,R)\cdot(B_1,\cdots,B_m)=
(\diag(C_1,0), \diag(C_2,0), \diag(C_3,0), 0, \cdots, 0,$$ where for some $\lambda\geq0$,
$$C_1:=
\begin{pmatrix}
0 & \lambda & 0 \\
-\lambda & 0 & 0 \\
0 & 0 & 0 \\
\end{pmatrix}, \
C_2:=
\begin{pmatrix}
0 & 0 & \lambda \\
0 & 0 & 0 \\
-\lambda & 0 & 0 \\
\end{pmatrix},\
C_3:=
\begin{pmatrix}
0 & 0 & 0 \\
0 & 0 & \lambda \\
0 & -\lambda & 0 \\
\end{pmatrix}.$$

\item If $n\geq 4$, then
\begin{equation*}
\sum^m_{r,s=1}\|\[B_r,B_s\]\|^2\leq \frac23\(\sum^m_{r=1}\|B_r\|^2\)^2.
\end{equation*}
The equality holds if and only if there exists a $(P,R)\in K(n, m)$ such that
$$(P,R)\cdot(B_1,\cdots,B_m)=
(\diag(D_1,0), \diag(D_2,0), \diag(D_3,0), 0, \cdots, 0),$$ where for some $\lambda\geq0$,
$$D_1:=
\begin{pmatrix}
0 & \lambda & 0 & 0 \\
-\lambda & 0 & 0 & 0 \\
0 & 0 & 0 & \lambda \\
0 & 0 & -\lambda & 0 \\
\end{pmatrix}, \
D_2:=
\begin{pmatrix}
0 & 0 & \lambda &0 \\
0 & 0 & 0 & -\lambda \\
-\lambda & 0 & 0 & 0 \\
0 & \lambda & 0 & 0 \\
\end{pmatrix},\
D_3:=
\begin{pmatrix}
0 & 0 & 0 & \lambda \\
0 & 0 & \lambda & 0 \\
0 & -\lambda & 0 & 0 \\
-\lambda & 0 & 0 & 0 \\
\end{pmatrix}.$$
\end{enumerate}
\end{thm}

For every $x\in M$, we denote by $\check{\kappa}(x)$ the largest eigenvalue of the curvature operator of $B$ at $\pi(x)\in B$, $\check{\lambda}(x)$ the lowest eigenvalue of the Ricci curvature of
$B$ at $\pi(x)\in B$ (thus $\check{\kappa}$ and $\check{\lambda}$ are constant along each fibre), and $\hat{\mu}(x)$ the largest eigenvalue of the Ricci curvature of the fibre at $x$.
Then we are ready to state the Simons-type integral inequalities derived by  the DDVV-type inequality for real skew-symmetric matrices.

\begin{thm}[Ge\cite{Ge14}]\label{geoskewDDVV}
Let $\pi: M^{n+m}\rightarrow B^n$ be a Riemannian submersion with totally geodesic fibres and
Yang-Mills horizontal distribution, and suppose $M$ is compact.
\begin{enumerate}
\item If $n=2$, then
\begin{equation*}
\int_M|A|^2\hat{\mu} dV_M\geq0.
\end{equation*}
\item If $m=1$, then
\begin{equation*}
\int_M|A|^2\left(\check{\kappa}-\check{\lambda}\right) dV_M\geq0.
\end{equation*}
\item If $m\geq 2$ and $n=3$, then
\begin{equation*}
\int_M|A|^2\left(\frac16|A|^2+2\hat{\mu}+\check{\kappa}-\check{\lambda}\right) dV_M\geq0.
\end{equation*}
\item If $m\geq 2$ and $n\geq 4$, then
\begin{equation*}
\int_M|A|^2\left(\frac13|A|^2+2\hat{\mu}+\check{\kappa}-\check{\lambda}\right) dV_M\geq0.
\end{equation*}
\end{enumerate}
\end{thm}

In \cite{Ge14}, the equality conditions for the above integral inequalities were characterized clearly.
Here we omit these equality conditions due to the complicated discussion.


\section{New progress on DDVV-type inequalities}
To describe the equality condition in this section, we put $\widetilde{K}(n, m):=U(n)\times O(m)$.
A $\widetilde{K}(n, m)$ action on a family of matrices $(A_1,\cdots,A_m)$ is given by
$$(P,R)\cdot(B_1,\cdots,B_m):=(P^*B_1P,\cdots,P^*B_mP)\cdot R.$$
\subsection{DDVV-type inequality for complex matrices}
In order to generalize the DDVV-type inequality to complex matrices,
Ge-Xu-You-Zhou \cite{GXYZ17} first considered the complex matrices with symmetries, namely, the Hermitian matrices and the skew-Hermitian matrices.

\begin{thm}[Ge-Xu-You-Zhou \cite{GXYZ17}]\label{Thm-Hermitian}
Let $B_1,\cdots,B_m$ be $n\times n$ (skew-)Hermitian matrices, $n\geq2$.
\begin{enumerate}[\rm (1)]
\item If\ $m\geq3$, then we have
$$\sum^m_{r,s=1}\|\[B_r,B_s\]\|^2\leq \frac43\(\sum^m_{r=1}\|B_r\|^2\)^2,$$
where the equality holds if and only if there exists a $(P,R)\in \widetilde{K}(n, m)$ such that
$$(P,R)\cdot(B_1,\cdots,B_m)=
(\diag(H_1,0), \diag(H_2,0), \diag(H_3,0), 0, \cdots, 0),$$ where for some $\lambda\geq0$,
$$H_1:=
\begin{pmatrix}
\lambda & 0 \\
0 & -\lambda \\
\end{pmatrix},
H_2:=
\begin{pmatrix}
0 & \lambda \\
\lambda & 0 \\
\end{pmatrix},
H_3:=
\begin{pmatrix}
0 & -\lambda\mathbf{i} \\
\lambda\mathbf{i} & 0 \\
\end{pmatrix}.$$

\item If $m=2$, then we have
$$\sum^2_{r,s=1}\|\[B_r,B_s\]\|^2\leq \(\sum^2_{r=1}\|B_r\|^2\)^2,$$
where the equality holds if and only if
under some $\widetilde{K}(n,2)$ action,
$B_1=diag(H_1,0)$ and $B_2=diag(\cos \theta H_2+\sin \theta H_3,0)$.
\end{enumerate}
\end{thm}

Using the  technique by dividing complex matrices into Hermitian matrices and skew-Hermitian matrices, Ge-Li-Zhou \cite{GLZ21} obtained the DDVV-type inequality for general complex matrices.
\begin{thm}[Ge-Li-Zhou \cite{GLZ21}]\label{TcDDVV}
Let $B_1,\cdots,B_m$ be arbitrary $n\times n$ complex matrices, $n\geq2$.
\begin{enumerate}[\rm (1)]
\item If $m\geq3$, then
\begin{equation*}\label{cDDVV}
\sum^m_{r,s=1}\|\[B_r,B_s\]\|^2\leq \frac43\(\sum^m_{r=1}\|B_r\|^2\)^2.
\end{equation*}
For $1\leq r \leq m$, let $B_r^1=\frac12(B_r+B_r^*), B_r^2=\frac12(B_r-B_r^*)$.
The equality holds if and only if $\sum^m_{r=1}\[B_r,B_r^*\]=0$ and
there exists a $(P,R)\in \widetilde{K}(n, 2m)$ such that
$$(P,R)\cdot(B_1^1, \cdots, B_m^1, \mathbf{i}B_1^2, \cdots, \mathbf{i}B_m^2)=
(\diag(H_1,0), \diag(H_2,0), \diag(H_3,0), 0, \cdots, 0),$$ where for some $\lambda\geq0$,
$$H_1:=
\begin{pmatrix}
\lambda & 0 \\
0 & -\lambda \\
\end{pmatrix}, \quad
H_2:=
\begin{pmatrix}
0 & \lambda \\
\lambda & 0 \\
\end{pmatrix},\quad
H_3:=
\begin{pmatrix}
0 & -\lambda\mathbf{i} \\
\lambda\mathbf{i} & 0 \\
\end{pmatrix}.
$$

\item If $m=2$, then
$$\sum^2_{r,s=1}\|\[B_r,B_s\]\|^2\leq \(\sum^2_{r=1}\|B_r\|^2\)^2.$$
The equality holds if and only if there exists a unitary matrix $U$
such that $B_1=U^*\diag(\widetilde{B_1}, 0)U, B_2=U^*\diag(\widetilde{B_2}, 0)U$,
where $\widetilde{B_1}, \widetilde{B_2}\in M(2, \mathbb{C})$
with $\|\widetilde{B_1}\|=\|\widetilde{B_2}\|,
\<\widetilde{B_1}, \widetilde{B_2}\>=0, \operatorname{tr}\widetilde{B_1}=\operatorname{tr}\widetilde{B_2}=0$.
\end{enumerate}
\end{thm}

\begin{rem}
Let
$$B_1:=
\begin{pmatrix}
1 & 0 \\
0 & -1 \\
\end{pmatrix}, \quad
B_2:=
\begin{pmatrix}
0 & 1 \\
1 & 0 \\
\end{pmatrix}, \quad
B_3:=
\begin{pmatrix}
0 & -1 \\
1 & 0 \\
\end{pmatrix},$$
then$$\sum^3_{r,s=1}\|\[B_r,B_s\]\|^2\ = \frac43\(\sum^3_{r=1}\|B_r\|^2\)^2, \quad \sum^2_{r,s=1}\|\[B_r,B_s\]\|^2\ = \(\sum^2_{r=1}\|B_r\|^2\)^2.$$
Hence the optimal constants for the real matrices case and the complex matrices case are both
$\frac43$ for $m\geq3$, and $1$ for $m=2$.
\end{rem}

By slightly changing the proof of Theorem \ref{TcDDVV},
Ge-Li-Zhou \cite{GLZ21} also obtained the DDVV-type inequality for complex symmetric or complex skew-symmetric matrices.
See Table \ref{table1} for the optimal constants.

\subsection{DDVV-type inequality for quaternionic matrices}\label{sec3}
Ge-Li-Zhou \cite{GLZ21} generalized the BW inequality and the DDVV inequality to quaternionic matrices.  
In this case, it turns out that both of
the optimal constants $c$ are double of that for complex matrices,
mainly because the multiplication of $\mathbf{i}, \mathbf{j}, \mathbf{k}$ is anti-commutative.
The proof is carried out simply by mapping a quaternionic matrix into a complex matrix  and then using the known inequalities for complex matrices
in a fashion analogous to what they have done with proving Theorem \ref{Thm-Hermitian}.

\begin{thm}[BW-type inequality for quaternionic matrices \cite{GLZ21}]\label{TqBW2}
Let $X,Y\in M(n,\mathbb{H})$, then $\|[X,Y]\|^2\leq 4\|X\|^2 \|Y\|^2$. The equality holds if and only if
either $\|X\| \|Y\|=0$, or $X=upu^*$ and $Y=uqu^*$ for some unit column vector $u\in\mathbb{H}^n$, where $p,q\in\mathbb{H}$ are 
purely imaginary quaternions that have real-orthogonal vector representations in the canonical basis.
\end{thm}

\begin{rem}
The maximal pair $(X,Y)$ can be rewritten as $X=U\diag(p,0,\cdots,0)U^*$, $Y=U\diag(q,0,\cdots,0)U^*$ for some quaternionic unitary matrix $U\in Sp(n)$. The condition on $p,q$ is equivalent to the anti-commutativity $pq=-qp$, which cannot happen in the real or complex cases.
\end{rem}

\begin{rem}
Let $X=\mathbf{i},~ Y=\mathbf{j}$, then
$$\|\[X,Y\]\|^2= 4\|X\|^2\|Y\|^2. $$
Hence $4$ is the optimal constant for the BW-type inequality for quaternionic matrices.
Moreover, for any $x\in\mathbb{R}^n$ and $\lambda\in\mathbb{R}$, let $X=xx^t\mathbf{i},~Y=\lambda xx^t\mathbf{j}\in M(n,\mathbb{H})$, we have
$$\|\[X,Y\]\|^2= 4\|X\|^2\|Y\|^2. $$
\end{rem}
In order to characterize the equality condition of the DDVV-type inequality for quaternionic matrices, we introduce the following map.
Let
$$\begin{aligned}
\Psi: M(n,\mathbb{H})&\longrightarrow M(2n,\mathbb{C}),\\
X=X_1+X_2\mathbf{j}&\longmapsto
\begin{pmatrix}
X_1 & X_2 \\
-\overline{X_2} & \overline{X_1} \\
\end{pmatrix},
\end{aligned}$$
where  $X_1, X_2 \in M(n,\mathbb{C})$.
It is easy to see
\begin{equation*}\label{qDDVV2}
\|\Psi(X)\|^2=2\|X\|^2.
\end{equation*}
For $X, Y \in M(n,\mathbb{H})$, let
$$A_1:=[X_1, Y_1]-X_2\overline{Y_2}+Y_2\overline{X_2}, \quad
A_2:=X_1Y_2-Y_2\overline{X_1}+X_2\overline{Y_1}-Y_1X_2.$$
Then direct calculations show
 $$[X,Y]=A_1+A_2\mathbf{j},\quad
\[\Psi(X),\Psi(Y)\]=
\begin{pmatrix}
A_1 & A_2 \\
-\overline{A_2} & \overline{A_1} \\
\end{pmatrix}. $$
Remember $\mathbf{j}\mathbf{i}=-\mathbf{i}\mathbf{j}$ when moving the $\mathbf{j}$ to the right.
Therefore $\Psi$ preserves the commutator:
\begin{equation*}\label{qDDVV3}
\[\Psi(X),\Psi(Y)\]=\Psi(\[X,Y\]).
\end{equation*}
\begin{thm}[Ge-Li-Zhou \cite{GLZ21}]\label{TqDDVV}
Let $B_1,\cdots,B_m$ be arbitrary $n\times n$ quaternionic matrices.
\begin{enumerate}[\rm (1)]
\item If $m\geq3$, then
\begin{equation*}\label{qDDVV}
\sum^m_{r,s=1}\|\[B_r,B_s\]\|^2\leq \frac83\(\sum^m_{r=1}\|B_r\|^2\)^2.
\end{equation*}
For $1\leq r \leq m$, let $B_r^1=\frac12(\Psi(B_r)+\Psi(B_r)^*), B_r^2=\frac12(\Psi(B_r)-\Psi(B_r)^*)$.
The equality holds if and only if $\sum^m_{r=1}\[\Psi(B_r),\Psi(B_r)^*\]=0$
and there exists a $(P,R)\in\widetilde{K}(n, 2m)$ such that
$$(P,R)\cdot(B_1^1, \cdots, B_m^1, \mathbf{i}B_1^2, \cdots, \mathbf{i}B_m^2)=
(\diag(H_1,0), \diag(H_2,0), \diag(H_3,0), 0, \cdots, 0),$$ where for some $\lambda\geq0$,
$$H_1:=
\begin{pmatrix}
\lambda & 0 \\
0 & -\lambda \\
\end{pmatrix}, \quad
H_2:=
\begin{pmatrix}
0 & \lambda \\
\lambda & 0 \\
\end{pmatrix}, \quad
H_3:=
\begin{pmatrix}
0 & -\lambda\mathbf{i} \\
\lambda\mathbf{i} & 0 \\
\end{pmatrix}.
$$

\item If $m=2$, then
$$\sum^2_{r,s=1}\|\[B_r,B_s\]\|^2\leq 2\(\sum^2_{r=1}\|B_r\|^2\)^2. $$
The equality holds if and only if  $B_1=upu^*$ and $B_2=uqu^*$ for some unit column vector $u\in\mathbb{H}^n$, where $p,q\in\mathbb{H}$ are orthogonal imaginary quaternions and $\|p\|=\|q\|$.
\end{enumerate}
\end{thm}

\begin{rem}
Let $B_1=\mathbf{i}, B_2=\mathbf{j}, B_3=\mathbf{k}$, then
$$\sum^3_{r,s=1}\|\[B_r,B_s\]\|^2\ = \frac83\(\sum^3_{r=1}\|B_r\|^2\)^2, $$
$$\sum^2_{r,s=1}\|\[B_r,B_s\]\|^2\ = 2\(\sum^2_{r=1}\|B_r\|^2\)^2. $$
Hence the optimal constants for the quaternionic matrices case
and the quaternionic skew-Hermitian matrices case are both
$\frac83$ for $m\geq3$, and $2$ for $m=2$. The equality condition could be also written in quaternion domain as in Theorem \ref{TqBW2}. A maximal triple $(B_1,B_2,B_3)$ determines the $(P, R)$-action, and all others must be zero. The surviving matrices should be in the form
$$B_r=uq_ru^*\in M(n,\mathbb{H}), \quad r=1,2,3,$$
for some unit column vector $u\in\mathbb{H}^n$ and $q_1,q_2,q_3\in\mathbb{H}$ are orthogonal imaginary quaternions with the same norm.
Also note that one can try to tackle the transition from the real to complex matrices in a similar, natural way. However, the doubled constant turns out not to be sharp in this case.
\end{rem}

\subsection{DDVV-type inequality for Clifford system and Clifford algebra}\label{sec4}
Inspired by the relation between DDVV-type inequalities and Erd\H{o}s-Mordell inequality (cf. \cite{E35, MB37, Ba58, Le61, DNP16, MM17,GS05}, etc) discovered by Z. Lu,
Ge-Li-Zhou \cite{GLZ21} established the DDVV-type inequalities for matrices in the subspaces spanned by a Clifford system or a Clifford algebra.

To illustrate the number $c$ more explicitly, we briefly introduce the representation theory of Clifford algebra (cf. \cite{FKM81}). A Clifford system on $\mathbb{R}^{2l}$ can be represented by real symmetric orthogonal matrices $P_0,\cdots,P_m\in O(2l)$ satisfying $P_iP_j+P_jP_i=2\delta_{ij}I_{2l}$; a Clifford algebra on $\mathbb{R}^l$ can be represented by real skew-symmetric orthogonal matrices $E_1,\cdots,E_{m-1}\in O(l)$ satisfying $E_iE_j+E_jE_i=-2\delta_{ij}I_{l}$; they are one-to-one correspondent by setting
\begin{equation*}
P_0=\begin{pmatrix}
    I_l &  0 \\
     0 & -I_l
\end{pmatrix}, \quad
P_1=\begin{pmatrix}
  0 & I_l    \\
 I_l & 0
\end{pmatrix}, \quad
P_{\alpha+1}=\begin{pmatrix}
  0 & E_\alpha    \\
 -E_\alpha & 0
\end{pmatrix}, \quad \alpha=1,\cdots,m-1.
\end{equation*}
A Clifford system $(P_0,\cdots,P_m)$ on $\mathbb{R}^{2l}$ (resp. Clifford algebra $(E_1,\cdots,E_{m-1})$  on $\mathbb{R}^l$) can be decomposed into a direct sum of $k$ irreducible Clifford systems (resp. Clifford algebras) on $\mathbb{R}^{2\delta(m)}$ (resp. on $\mathbb{R}^{\delta(m)}$) with $l=k\delta(m)$ for $k,m\in\mathbb{N}$, where the irreducible dimension $\delta(m)$ satisfies $\delta(m+8)=16\delta(m)$ and can be listed in the following Table \ref{table3}.
\begin{table}[h!]
\caption{Dimension $\delta(m)$ of irreducible representation of Clifford algebra}\label{table3}
\centering
\begin{tabular}{|c|c|c|c|c|c|c|c|c|c|}
\hline
$m$ & $1$ & $2$ & $3$ & $4$ & $5$ & $6$ & $7$ & $8$ & $\cdots~ m+8$\\
\hline
$\delta(m)$& $1$ &$2$ &$4$ &$4$ & $8$ & $8$ & $8$ & $8$ & $\cdots~ 16\delta(m)$\\
\hline
\end{tabular}
\end{table}

\begin{thm}[Ge-Li-Zhou \cite{GLZ21}]\label{TCsDDVV}
Let $(P_0, P_1, \cdots, P_m)$ be a Clifford system on $\mathbb{R}^{2l}$, i.e., $P_0,\cdots,P_m\in O(2l)$ are real symmetric orthogonal matrices satisfying $P_iP_j+P_jP_i=2\delta_{ij}I_{2l}$.
Let $B_1,\cdots,B_M \in span\{P_0, P_1, \cdots, P_m\}$, then
\begin{equation*}\label{CsDDVV}
\sum^M_{r,s=1}\|\[B_r,B_s\]\|^2\leq \frac2l\(1-\frac1N\)\(\sum^M_{r=1}\|B_r\|^2\)^2, \quad
N=\min\{m+1, M\}.
\end{equation*}
The condition for equality has two cases.

Case $(1)$: When $m+1\leq M$, the equality holds if and only if
$p_0, \cdots, p_{m}$ are orthogonal vectors with the same norm,
where $p_i:=\(\<P_i, B_1\>, \cdots, \<P_i, B_M\>\)\in\mathbb{R}^M$.

Case $(2)$: When $m+1\geq M$, the equality holds if and only if
$B_1,\cdots,B_M$ are orthogonal matrices with the same Frobenius norm.

\end{thm}

Analogously they were able to obtain the DDVV-type inequality for Clifford algebra.
\begin{thm}[Ge-Li-Zhou \cite{GLZ21}]\label{TCaDDVV}
Let $\{E_1, \cdots, E_{m-1}\}$ be a Clifford algebra on $\mathbb{R}^l$, i.e.,  $E_1,\cdots,E_{m-1}\in O(l)$ are real skew-symmetric orthogonal matrices satisfying $E_iE_j+E_jE_i=-2\delta_{ij}I_{l}$. Let $B_1,\cdots,B_M \in span\{E_1, \cdots, E_{m-1}\}$, then
\begin{equation*}\label{CaDDVV}
\sum^M_{r,s=1}\|\[B_r,B_s\]\|^2\leq \frac4l\(1-\frac1N\)\(\sum^M_{r=1}\|B_r\|^2\)^2,\quad
N=\min\{m-1, M\}.
\end{equation*}
The condition for equality has two cases.

Case $(1)$: When $m-1\leq M$, the equality holds if and only if
$e_1, \cdots, e_{m-1}$ are orthogonal with the same norm,
where $e_i:=\(\<E_i, B_1\>, \cdots, \<E_i, B_M\>\)\in\mathbb{R}^M$.

Case $(2)$: When $m-1\geq M$, the equality holds if and only if
$B_1,\cdots,B_M$ are orthogonal with the same norm.

\end{thm}

\subsection{Lu inequality}
To prove DDVV conjecture, Lu \cite{Lu11} discovered the stronger inequality below. Meanwhile, this inequality can also {provide a new geometric pinching result (see Theorem \ref{Lupin}).}
\begin{thm}[{Lu inequality \cite{Lu11}}]\label{thm Lu inequality}
Let $A$ be an $n\times n$  diagonal  matrix of norm $1$. Let $A_2,\cdots, A_m$ be symmetric matrices such that
\begin{enumerate}[\rm (i)]
\item $\langle A_\alpha,A_\beta\rangle=0$ if  $\alpha\neq \beta$;
\item $\|A_2\|\geq\cdots\geq\|A_m\|$.
\end{enumerate}
Then we have
\begin{equation}\label{key-1}
\sum_{\alpha=2}^m\|[A,A_\alpha]\|^2\leq \sum_{\alpha=2}^m \|A_\alpha\|^2+\|A_2\|^2.
\end{equation}
The  equality in~\eqref{key-1} holds if and only if, after an orthonormal base change and up to a sign, we have
\begin{enumerate}[\rm (1)]
\item
$A_3=\cdots =A_m=0$, and
\begin{equation*}\label{a116}
A=\begin{pmatrix}
\frac{1}{\sqrt 2} &0\\0&-\frac{1}{\sqrt 2}\\
&&0\\
&&&\ddots\\
&&&&0
\end{pmatrix},\quad
A_2=c\begin{pmatrix}
0&\frac{1}{\sqrt 2}\\\frac{1}{\sqrt 2}&0\\
&&0\\
&&&\ddots\\
&&&&0
\end{pmatrix},
\end{equation*}
where $c$ is any constant, or
\item For two real numbers $\lambda=1/\sqrt{n(n-1)}$ and $\mu$, we have
\begin{equation*}\label{a112}
A=\lambda\begin{pmatrix}
n-1\\
&-1\\
&&\ddots\\
&&&-1
\end{pmatrix},
\end{equation*}
and
$A_\alpha$ is $\mu$ times the matrix whose only nonzero entries are $1$ at the $(1,\alpha)$ and $(\alpha,1)$ places, where $\alpha=2,\cdots,n$.
\end{enumerate}
\end{thm}
 Equivalently, one has the following geometric version of the Lu inequality.
 \begin{thm}[Normal Ricci curvature inequality \cite{Lu11}]
  For any $1\leq r \leq m$,
 \begin{equation*}
 {\rm Ric}^{\bot}\left(\xi_r, \xi_r\right)
 \leq \(\max \limits_{s\not=r}\|B_s\|^2+\sum_{s\not=r}\|B_s\|^2\)\|B_r\|^2,
 \end{equation*}
 where $ {\rm Ric}^{\bot}\left(\xi_r, \xi_r\right):=\left(\sum_{s\not=r}\| R^{\bot}_{rs}\|^2 \right) ^{\frac{1}{2}}$ is the  normal Ricci curvature.
 \end{thm}

\section{Open questions}
In this section, we would like to mention more about possible future studies on DDVV-type inequalities and related topics.
\subsection{Questions about DDVV-type inequalities}
We begin with questions about the algebraic aspects of DDVV-type inequalities.
\begin{ques}
What can we expect for the quotient
$$\frac{\sum_{r,s=1}^m\|[B_r,B_s]\|^2}{\left(\sum_{r=1}^{m}\|B_r\|^2 \right) ^2}?$$
What is the expectation of the commutators of random matrices in certain categories like GOE, GUE, and GSE?
\end{ques}
It seems that the lists of the optimal constant $c$ would possibly have some links with random matrix theory or quantum physics. However, it is just our naive and wild guess since we know nothing about that.

\begin{ques}
Find more DDVV-type inequalities for matrices, Lie algebras or operators lying in certain subspaces of special interest like spaces of austere matrices (see for a special example in \cite{GTY20}).
\end{ques}

L\'{a}szl\'{o} \cite{Laszlo2010} proved that: \emph{the smallest value $\gamma_n$ so that the nonnegative polynomial}
\begin{equation*}\label{equation Lu-Wenzel-type SOS}
F(X,Y):=\left( 2+\gamma_n\right) \left( \|X\|^2\|Y\|^2-|\left\langle X,Y \right\rangle |^2\right) -\|[X,Y]\|^2 \quad\quad  (X, Y \in \mathbb{R}^{n\times n})
\end{equation*}
\emph{is a sum of squares (SOS) of polynomials is $\gamma_n=\frac{n-2}{2}$.}
Based on this result and many similar results (cf.  \cite{Laszlo2012,LW17}), Lu-Wenzel \cite{LW17} proposed the following conjecture.
\begin{conj}[Lu-Wenzel \cite{LW17}]
The form
$$2\|X\|^2\|Y\|^2-2|\left\langle X,Y \right\rangle |^2-\|[X,Y]\|^2$$
generated by two arbitrary real Toeplitz matrices is SOS.
\end{conj}
We are also concerned about the analogous question for DDVV-type inequalities:
\begin{ques}
Whether the nonnegative polynomial defined by the DDVV-type inequalities (by $F(B_1,\cdots,B_m):=c\left(\sum_{r=1}^{m}\|B_r\|^2 \right) ^2-\sum_{r,s=1}^m\|[B_r,B_s]\|^2$)
is a sum of squares of quadratic forms on the matrices in the regarded types?  
This would provide more examples on the generalized Hilbert's 17th problem (cf. \cite{GT22}).
\end{ques}

The next two questions are related to the geometric aspect of the DDVV inequality.
\begin{ques}
What can we expect for a minimal submanifold in $\mathbf{S}^{n+1}$ with normal scalar curvature pinched?
\end{ques}
In \cite{GTY20}, Ge-Tang-Yan obtained new normal scalar curvature inequalities
(which is sharper than the DDVV inequality) on the focal submanifolds of isoparametric hypersurfaces in the unit sphere, and they characterized the subsets which achieve upper or lower bounds.

\begin{ques}
Classify all the submanifolds that the DDVV inequality achieves equality at every point.
\end{ques}
Submanifolds achieving equality of the DDVV inequality everywhere are called Wintgen ideal submanifolds, which are not classified so far (cf. \cite{CL08, DT09, LMW15, LMWX16, WX14, X15, X17, XLMW14, XLMW18}).
See \cite{C21} for a detailed survey on the research of Wintgen ideal submanifolds.

\subsection{Generalized Peng-Terng 2nd-Gap.}
By Theorem \ref{Lupin},
the quantity $S+\lambda_2$ might be the right object to study pinching theorems. To justify this, we introduce the following Lu's conjecture:

\begin{conj}[Lu \cite{Lu11}]\label{pt}
Let $M$ be an $n$-dimensional closed minimal submanifold in the unit sphere $\mathbf{S}^{n+m}$. If $S+\lambda_2$ is a constant and if
$$S+\lambda_2>n,$$
then there is a constant $\varepsilon(n,m)>0$ such that
$$S+\lambda_2>n+\varepsilon(n,m).$$
\end{conj}

If $m=1$, this conjecture is true (cf. \cite{PT83, DX11, LXX21}). In fact, this is a special case of Chern's conjecture (cf. \cite{Chern68, CdK}). For more details, please refer to  \cite{Chang93, PT83-2, TY20, SWY12}, etc.

\subsection{Lu-Wenzel Conjectures}
In order to give a unified generalization of the BW inequality, DDVV inequality and Lu  inequality,
Lu and Wenzel  (\cite{LW16, LW17})  proposed several conjectures (also called LW Conjectures) in 2016.
They started with the following 
 Conjectures \ref{conj1}, \ref{conj2}, \ref{conj3} and an open Question \ref{Lu-Wenzel ques} in the space $M(n, \mathbb{K})$
of $n\times n$ matrices over the field $\mathbb{K}=\mathbb{R}$.
\begin{conj}[Lu-Wenzel \cite{LW16, LW17}]
\label{conj1}
Let $B_1,\cdots,B_m\in M(n,\mathbb{K})$ 
subject to  $$\operatorname{tr}\Big(B_\alpha [B_\gamma, B_\beta]\Big)=0$$ for any $1\leq \alpha,\beta,\gamma\leq m$, then
\begin{equation*}\label{csymDDVV}
\sum^m_{\alpha,\beta=1}\|\[B_\alpha,B_\beta\]\|^2\leq  \(\sum^m_{\alpha=1}\|B_\alpha\|^2\)^2.
\end{equation*}
\end{conj}

\begin{conj}[Fundamental Conjecture of Lu-Wenzel \cite{LW16, LW17}]\label{conj2}
Let $B,B_2,\cdots,B_m\in M(n,\mathbb{K})$ be matrices such that
\begin{enumerate} [\rm (i)]
\item $\operatorname{tr}(B_\alpha B^*_\beta)=0$,
i.e., $B_\alpha\bot B_\beta$ for any $\alpha\neq\beta$;
\item $\operatorname{tr}\Big(B_\alpha [B, B_\beta]\Big)=0$ for any $2\leq \alpha,\beta\leq m$.
\end{enumerate}
 Then
\begin{equation*}\label{ineq of conj2}
\sum^m_{\alpha=2}\|\[B,B_\alpha\]\|^2\leq \(\max \limits_{2\leq \alpha\leq m}\|B_\alpha\|^2+\sum^m_{\alpha=2}\|B_\alpha\|^2\)\|B\|^2.
\end{equation*}
\end{conj}

Note that the Lu inequality is a  special case of Conjecture \ref{conj2}. 

Let's consider the linear operator $T_X$ as in the following conjectures and questions.  
 More specifically, for any $n\times n$ complex matrix $X$ with $\|X\|=1$, we define
$$\begin{aligned}
T_X:  M(n,\mathbb{C})&\longrightarrow M(n,\mathbb{C}),\\
Y&\longmapsto [X^*,[X,Y]].
\end{aligned}$$
It turns out that $T_X$ is exactly an operator on $V=M(n,\mathbb{C})$ and  $\dim_{\mathbb{C}} V=n^2$. 

\begin{prop}[\cite{GLLZ20}]\label{prop0}
$T_X$ has the following properties:
\begin{enumerate}
\item[$(a)$] $T_X$ is a self-dual and positive semi-definite linear map.
\item[$(b)$] The set of eigenvalues $\lambda(T_X):=\{\lambda_1(T_X)\geq\cdots \geq\lambda_N(T_X)\}$ is invariant under unitary congruences of $X$.
\item[$(c)$] The multiplicity of each positive eigenvalue of $T_X$ is even, i.e., $\lambda_{2i-1}(T_X)=\lambda_{2i}(T_X)$ for any $i$ with $\lambda_{2i-1}(T_X)>0$.
\end{enumerate}
\end{prop}

\begin{conj}[Lu-Wenzel \cite{LW16, LW17}]
\label{conj3}
For $X\in M(n,\mathbb{K})$ with $\|X\|=1$,  
then
$$\lambda_1(T_X)+\lambda_3(T_X)\leq3.$$
\end{conj}

\begin{ques}[Lu-Wenzel \cite{LW16, LW17}]
\label{Lu-Wenzel ques}
What is the upper bound of
$\sum\limits_{i=1}^{k}\lambda_{2i-1}(T_X)?$
\end{ques}
\begin{rem}
In Question \ref{Lu-Wenzel ques}, one has
\begin{enumerate}
\item
If $k = 1$, the bound is $2$ by the BW inequality, i.e., $\lambda_1(T_X)\leq2$, since 
$$\lambda_1(T_X)=\max \limits_{\|Y\|=1}\langle T_X Y, Y\rangle=\max \limits_{\|Y\|=1}\|[X,Y]\|^2\leq 2.$$
\item If $k = 2$, the bound is supposed to be $3$ by Conjecture \ref{conj3}.
\end{enumerate}
\end{rem}
How are all these conjectures and the known inequalities connected?
When restricted to real symmetric matrices, Conjecture \ref{conj1} reduces to the DDVV inequality. It turns out that not only the BW inequality and the DDVV inequality but also both Conjectures \ref{conj1} and \ref{conj3} are implied by Conjecture \ref{conj2} (cf. \cite{LW16}).
Hence, Conjecture \ref{conj2} (as well as the equivalent Conjectures \ref{conj2A}--\ref{conj2C}, see  Theorem \ref{equiv conj}) takes exactly the role of a unified generalization of the BW inequality and the DDVV inequality for real matrices. We call it the Fundamental Conjecture of Lu and Wenzel, or simply the  (real, i.e., $\mathbb{K}=\mathbb{R}$)
LW Conjecture. 

Next, we  will introduce some equivalent forms of Conjecture \ref{conj2}.
Since the BW inequality (resp.  the DDVV inequality) holds also for complex (resp. complex symmetric) matrices (cf. \cite{BW08}, \cite{GLZ21}), we can also consider the same conjectures as above in the complex version. 
Let $\mathbb{K}=\mathbb{R}$ or $\mathbb{K}=\mathbb{C}$ in the following Conjectures \ref{conj2A} - \ref{conj2C}.
\begin{conj}[Ge-Li-Lu-Zhou \cite{GLLZ20}]\label{conj2A}
For $X\in M(n,\mathbb{K})$ with $\|X\|=1$, we have
\begin{equation*}\label{ineq of conj2A}
\sum_{i=1}^{2k}\lambda_{i}(T_X)\leq2k+2, \quad k=1,\cdots,\[\frac{n^2}{2}\].
\end{equation*}
\end{conj}

In fact, the sum $\sum_{i=1}^{2k}\lambda_{i}(T_X)$ in Conjecture \ref{conj2A} cannot exceed $2n$. We explain this by introducing the following Conjecture \ref{conj2B} which looks stronger but in fact is equivalent to Conjecture \ref{conj2A}.


\begin{conj}[Ge-Li-Lu-Zhou \cite{GLLZ20}]\label{conj2B}
For $X\in M(n,\mathbb{K})$ with $\|X\|=1$,
we have
\begin{align*}\label{ineq of conj2C}
\begin{split}
\sum_{i=1}^{2k}\lambda_{i}(T_X)\leq \left \{
\begin{array}{ll}
2k+2,       & 1 \leq k\leq n-1;\\
2n,       & n \leq k.\\
\end{array}
\right.
\end{split}
\end{align*}
\end{conj}
Another equivalent conjecture that also appears to be stronger is the following Conjecture \ref{conj2C}.
It omits the second assumption of Conjecture \ref{conj2}, at the price of a factor 2 in the bound.

\begin{conj}[Ge-Li-Lu-Zhou \cite{GLLZ20}]\label{conj2C}
Let $B,B_2,\cdots,B_m\in M(n,\mathbb{K})$ be matrices
such that $\operatorname{tr}(B_\alpha B^*_\beta)=0$ for any $2\leq\alpha\neq \beta\leq m$. Then
$$\sum^m_{\alpha=2}\|\[B,B_\alpha\]\|^2\leq \(2\max_{2\leq \alpha\leq m}\|B_\alpha\|^2+\sum^m_{\alpha=2}\|B_\alpha\|^2\)\|B\|^2.$$
\end{conj}

We summarize the relations of these conjectures in the following theorem.

\begin{thm}[Ge-Li-Lu-Zhou \cite{GLLZ20}]\label{equiv conj}
If  $\mathbb{K}=\mathbb{R}$ or $\mathbb{K}=\mathbb{C}$, 
then
the following relations hold in these conjectures.
\begin{enumerate} [ \rm (i) ]
\item Conjectures \ref{conj2}, \ref{conj2A}, \ref{conj2B}, and \ref{conj2C} are equivalent to each other.
\item If one of the  conjectures above is true, then Conjectures \ref{conj1} and \ref{conj3} hold.
\end{enumerate}
\end{thm}

%
%
%

Hence, we call Conjecture \ref{conj2} for complex matrices (i.e., $\mathbb{K}=\mathbb{C}$) the complex LW Conjecture. Obviously, the complex LW conjecture implies the real LW conjecture.
Ge-Li-Lu-Zhou \cite{GLLZ20} proved the complex LW Conjecture  in some special cases which we conclude in the following.
\begin{thm}[Ge-Li-Lu-Zhou \cite{GLLZ20}]\label{thm-special-cases}
The complex LW Conjectures \ref{conj2A} and \ref{conj2B}
(and because of Theorem \ref{equiv conj} all conjectures of this section)
are true in one of the following cases:
\begin{enumerate} [ \rm (i) ]
 \item $X\in M(n,\mathbb{C})$ is a normal matrix;
 \item $\rank X=1$;
 \item $n=2$ or $n=3$.
\end{enumerate}
\end{thm}

For the Conjectures \ref{conj2A} and \ref{conj2B}
in general Ge-Li-Lu-Zhou \cite{GLLZ20} were able to get some weakened results as follows.
\begin{thm}[Ge-Li-Lu-Zhou \cite{GLLZ20}]\label{thm-conj3-weak}
For $X\in M(n,\mathbb{C})$ with $\|X\|=1$, we have
\begin{enumerate} [ \rm (i) ]
 \item
 $\lambda_1(T_X)+\lambda_3(T_X)\leq\frac{4+\sqrt{10}}{2}\approx 3.58;$

 \item
 $\sum_{t=1}^{2k}\lambda_{i}(T_X)\leq2k+1+2\sqrt{k},\quad k=1,\cdots,\[\frac{n^2}{2}\].$
 \end{enumerate}
\end{thm}


\begin{ack}
The authors would like to thank the referee for valuable comments and suggestions.
\end{ack}


\end{document}